\documentclass[a4paper,11pt]{article}

\usepackage{amsfonts,amssymb,amsmath,amsthm,latexsym}
\usepackage{mathrsfs,bbm}
\usepackage{graphicx,graphics,psfrag,epsf,eepic}

\usepackage[english]{babel}
 \makeatletter
 \def\@captype{figure}
 \makeatother

 \makeatletter
 \@addtoreset{equation}{section}
 \makeatother

 \makeatletter
 \@addtoreset{enunciato}{section}
 \makeatother

 \newcounter{enunciato}[section]

 \newtheorem{ittheorem}{Theorem}
 \newtheorem{itlemma}{Lemma}
 \newtheorem{itproposition}{Proposition}
 \newtheorem{itdefinition}{Definition}
 \newtheorem{itremark}{Remark}
 \newtheorem{itclaim}{Claim}
 \newtheorem{itfact}{Fact}
 \newtheorem{itconjecture}{Conjecture}

 \newenvironment{theorem}{\addtocounter{enunciato}{1}
 \begin{ittheorem}}{\end{ittheorem}}

 \newenvironment{lemma}{\addtocounter{enunciato}{1}
 \begin{itlemma}}{\end{itlemma}}

 \newenvironment{proposition}{\addtocounter{enunciato}{1}
 \begin{itproposition}}{\end{itproposition}}

 \newenvironment{definition}{\addtocounter{enunciato}{1}
 \begin{itdefinition}}{\end{itdefinition}}

 \newenvironment{remark}{\addtocounter{enunciato}{1}
 \begin{itremark}}{\end{itremark}}

 \newenvironment{claim}{\addtocounter{enunciato}{1}
 \begin{itclaim}}{\end{itclaim}}

 \newenvironment{fact}{\addtocounter{enunciato}{1}
 \begin{itfact}}{\end{itfact}}

 \newenvironment{conjecture}{\addtocounter{enunciato}{1}
 \begin{itconjecture}}{\end{itconjecture}}

 \newcommand{\be}[1]{\begin{equation}\label{#1}}
 \newcommand{\ee}{\end{equation}}

 \newcommand{\bl}[1]{\begin{lemma}\label{#1}}
 \newcommand{\el}{\end{lemma}}

 \newcommand{\br}[1]{\begin{remark}\label{#1}}
 \newcommand{\er}{\end{remark}}

 \newcommand{\bt}[1]{\begin{theorem}\label{#1}}
 \newcommand{\et}{\end{theorem}}

 \newcommand{\bd}[1]{\begin{definition}\label{#1}}
 \newcommand{\ed}{\end{definition}}

 \newcommand{\bcl}[1]{\begin{claim}\label{#1}}
 \newcommand{\ecl}{\end{claim}}

 \newcommand{\bfact}[1]{\begin{fact}\label{#1}}
 \newcommand{\efact}{\end{fact}}

 \newcommand{\bp}[1]{\begin{proposition}\label{#1}}
 \newcommand{\ep}{\end{proposition}}

 \newcommand{\bc}[1]{\begin{corollary}\label{#1}}
 \newcommand{\ec}{\end{corollary}}

 \newcommand{\bcj}[1]{\begin{conjecture}\label{#1}}
 \newcommand{\ecj}{\end{conjecture}}

 \newcommand{\bpr}{\begin{proof}}
 \newcommand{\epr}{\end{proof}}

 \newcommand{\bprl}[1]{\begin{proofof}{\it\ref{#1}}.\,\,}
 \newcommand{\eprl}{\end{proofof}}

 \newcommand{\bi}{\begin{itemize}}
 \newcommand{\ei}{\end{itemize}}

 \newcommand{\ben}{\begin{enumerate}}
 \newcommand{\een}{\end{enumerate}}

 \newenvironment{proofof}{\noindent {\em Proof of Lemma\,\,}}{\hspace*{\fill}$\halmos$\medskip}
 \newcommand{\halmos}{\rule{1ex}{1.4ex}}

\def \cc {{\rm c}}

\def \Z {{\mathbb Z}}

\def \cE {{\mathcal E}}
\def \cF {{\mathcal F}}

\newcommand{\reff}[1]{(\ref{#1})}

\def\one{\rlap{\mbox{\small\rm 1}}\kern.15em 1}

\def\embf#1{\emph{\bf #1}}

\begin{document}
\title{Regular $g$-measures are not always Gibbsian}

\author{\renewcommand{\thefootnote}{\arabic{footnote}}
R.\ Fern\'andez
\footnotemark
\\
\renewcommand{\thefootnote}{\arabic{footnote}}
S.\ Gallo
\footnotemark
\\
\renewcommand{\thefootnote}{3,4}
G.\ Maillard
\footnotemark[4]
}

\footnotetext[1]{ 
Department of Mathematics, University of Utrecht,
P.O.box 80010, NL-3508 TA Utrecht, The Netherlands\\
{\sl R.Fernandez1@uu.nl}
}
\footnotetext[2]{
Instituto de Matem\'atica, Estat\'{\i}stica e Computa\c{c}\~ao Cient\'{\i}fica,
Universidade Estadual de Campinas,
Rua Sergio Buarque de Holanda, 651 13083-859 Campinas, Brasil\\
{\sl gsandro@ime.unicamp.br}
}
\footnotetext[3]{
CMI-LATP, Universit\'e de Provence,
39 rue F. Joliot-Curie, F-13453 Marseille Cedex 13, France\\
{\sl maillard@cmi.univ-mrs.fr}
}
\footnotetext[4]{
EURANDOM, P.O.\ Box 513, NL-5600 MB Eindhoven, The Netherlands
}
\date{}
\maketitle


\begin{abstract}
Regular $g$-measures are discrete-time processes determined by conditional 
expectations with respect to the past.  One-dimensional Gibbs measures, on the 
other hand, are fields determined by simultaneous conditioning on past and future.  
For the Markovian and exponentially continuous cases both theories are known 
to be equivalent.  Its equivalence for more general cases was an open problem.  
We present a simple example settling this issue in a negative way: there exist 
$g$-measures that are continuous and non-null but are not Gibbsian.  Our example 
belongs, in fact, to a well-studied family of processes with rather nice attributes:  It is 
a chain with variable-length memory, characterized by the absence of phase coexistence 
and the existence of a visible renewal scheme.   

{\it Acknowledgment.} GM is grateful to CNRS for financial support and to EURANDOM 
for hospitality. SG is grateful to EURANDOM and the Mathematics Department of the
University of Utrecht for hospitality. 
SG is supported by a FAPESP fellowship 2009/09809-1.  
This work is part of USP project \emph{Mathematics, computation, 
language and the brain}. We thank Aernout van Enter for his critical reading of the manuscript 
and his suggestions.
\end{abstract}



\section{Introduction}
\label{S1}

Measures on $E^{\mathbb Z}$ are the object of two very developed theories.  On the one hand, 
the theory of \emph{chains of complete connections}, started in~\cite{onimih35} and developed 
under a variety of names in slightly non-equivalent frameworks: \emph{chains of infinite 
order}~\cite{har55},  \emph{$g$-measures}~\cite{kea72}, \emph{uniform martingales}~\cite{kal90}, 
etc.  On the other hand, the theory of one-dimensional \emph{Gibbs measures} started 
in~\cite{dob68,lanrue69} and whose classical reference is the treatise~\cite{geo88}.  The former 
theory interprets $\mathbb Z$ as discrete time and measures as discrete-time processes.  The 
building blocks of the theory are, therefore, \emph{transition probabilities}, that is, conditional 
probabilities \emph{with respect to the past}. Chains, $g$-measures, etc are defined by their 
invariance under ---or consistency with--- these transition probabilities.  In contrast, Gibbsianness 
refers to fields in a spatial setting determined by distributions on finite-regions \emph{conditioned 
on exterior configurations}. In one dimension, if $\Z$ is interpreted as time, this corresponds to 
conditioning both \emph{with respect to the past and to the future}.  Of course, this is only part 
of the story.  Both processes and Gibbs measures are required to satisfy suitable regularity 
conditions, as reviewed below.

Given this state of affairs it is natural to wonder about the relation of both theories.  Are they 
equivalent? Does one-side conditioning carry the same information as two-side conditioning?  
Is every regular process a Gibbs measure and vice-versa? To be sure, there is another source 
of difference related to non-nullness.  The standard theory of Gibbs measures deals with systems 
with no forbidden configurations.  Many important instances of processes, on the other hand, 
include grammars or local exclusion rules (subshifts of finite type).  The previous questions should 
be stated, then, in the common non-null framework.  In this set-up, the equivalence of both 
theories has long been known to be true for Markov processes and fields (see, for instance, 
\cite[Chapter 11]{geo88}) and when continuity rates are exponentially decreasing~\cite{fermai04}.  
In this note we exhibit a simple example ---where all calculations can be explicitly performed--- 
showing that this equivalence is \emph{not} true in general.  

In more detail, below we construct a $g$-measure $\mu$ on $\{0,1\}^{\mathbb Z}$ with the 
following properties:
\begin{itemize}
\item $\mu$ is non-null: it gives nonzero measure to every cylinder.
\item There exists a continuous $g$ function for which $\mu$ is the \emph{only} consistent 
measure.
\item $\mu$ is left-right symmetric [as proven in the third line of \reff{ratioeq}].
\item $\mu$ is a variable-length memory chain which admits a renewal construction with visible renewals (see  \cite{gal11} for example).  
\item $\mu$ can be perfectly simulated (for example by the method of \cite{gal11}).
\end{itemize}
Yet, despite all these fine properties, the measure $\mu$ is non-Gibbsian. 

This example shows that regular $g$-measures are a different type of creature than Gibbsian 
measures and respective theories (large-deviations, uniqueness theorems, variational approach) 
can not, in general, be imported from one to the other.  In particular, readers are warned that a 
$g$-function of the form $g={\rm e}^\phi$ with $\phi$ ``nice'' does not automatically deserve the 
qualifier ``Gibbsian''.   These observations complement previous studies on  differences between one-sided and two-sided measurability done in the more general framework of ergodic theory
(see, e.g.~\cite{entver10} and references therein).

Another potential source of confusion arises from the traditional use, by people working in 
dynamical systems, of the word ``Gibbsian'' to refer to SRB measures (Sinai-Ruelle-Bowen 
measures, see e.g.\ \cite{rue78} or \cite{bow75}).  As briefly reviewed below, the set of non-null 
SRB measures is \emph{strictly contained} in the set of statistical mechanical Gibbs measures 
(but the former, unlike the latter, can also incorporate exclusions and subshifts).  Hence our 
$g$-measure $\mu$ is also non-Gibbsian in SRB sense.  


Let us conclude with a brief explanation of the non-Gibbsianness argument below.  The measure 
$\mu$ is non-null and consistent with a continuous $g$-measure.  This means that, upon 
conditioning, it becomes asymptotically insensitive to the far past.  In order to be Gibbs, the same 
asymptotic insensitivity must hold but simultaneously with respect to past \emph{and future}.  The 
measure $\mu$ is supported on configurations with an infinite number of $1$'s.   This is because 
the probability of having a $1$ after an infinite sequence of $0$'s is a strictly positive number 
$p_\infty$.   In addition, by continuity, the probability of having a $1$ conditioned on a large string 
of zeros converges, as the first $1$ recedes to $-\infty$, to $p_\infty$.  For the measure to be 
Gibbsian this same continuity must hold for \emph{two-side} conditioning.  A delicate case arises, 
however, when conditioning on having the all-zero configuration both towards the past and the 
future.  Of course, this is an impossible configuration for $\mu$, so that the actual value of this 
conditional probability is irrelevant.  Rather, the Gibbsianness question refers to  whether 
conditional probabilities converge (to whatever) as both the first $1$ to the left and the first $1$ 
to the right move away.  This is an essential property in the sense that its absence can not be fixed by 
measure-zero redefinitions.  Theorem \ref{th:1} shows that for some $g$-functions this two-side 
continuity is impossible.  As an aside, we point out that the ``all-0'' configuration is the \emph{only} point of discontinuity of the two-side conditional probabilities.  The measure $\mu$ is, therefore,  \emph{almost Gibbsian} and thus \emph{weakly Gibbsian} (see, for instance, \cite[Section 4.4]{fer06} for the corresponding definitions and historical references to these notions).


\section{Preliminaries}
\label{S2}

We consider a measurable space $(E,\cE)$ where $E=\{0,1\}$ is a two-symbol alphabet
and $\cE$ is the associated discrete $\sigma$-algebra.  We denote by $(\Omega,\cF)$ 
the associated product measurable space, that is $\Omega=E^{\Z}$, $\cF=\cE^{\Z}$.  
For each $\Lambda\subset\Z$ we denote $\Omega_\Lambda = E^{\Lambda}$ and 
$\sigma_\Lambda$ for the restriction of a configuration $\sigma\in\Omega$ to 
$\Omega_\Lambda$, namely the family $(\sigma_i)_{i\in\Lambda}\in E^{\Lambda}$. 
Also, $\cF_\Lambda$ will denote the sub-$\sigma$-algebra of $\mathcal{F}$ generated 
by cylinders based on $\Lambda$ ($\cF_\Lambda$-measurable functions are insensitive 
to configuration values outside $\Lambda$).  When $\Lambda$ is an interval, 
$\Lambda=[k,n]$ with $k,n\in\Z$ such that $k\le n$, we use the notation: $\omega_{k}^{n}
=\omega_{[k,n]}=\omega_{k},\ldots,\omega_{n}$, $\Omega_k^n=\Omega_{[k,n]}$ and 
$\cF_k^n=\cF_{[k,n]}$.  For semi-intervals we denote also $\cF_{\le n}=\cF_{(-\infty,n]}$, 
etc.  The concatenation notation $\omega_\Lambda\,\sigma_\Delta$, where 
$\Lambda\cap\Delta=\emptyset$, indicates the configuration on $\Lambda\cup\Delta$
coinciding with $\omega_i$ for $i\in\Lambda$ and with $\sigma_i$ for $i\in\Delta$.

We start by briefly reviewing the well-known notions of chains in a shift-invariant setting.
In this particular case, chains are also called $g$-measures (see \cite{kea72}).

\bd{lis1}
A \embf{regular $g$-function} $P$ on $\Omega$ is a probability kernel
$P\colon E\times\Omega_{-\infty}^{-1}\to [0,1]$, i.e.,
\be{g-funct}
\sum_{\omega_0\in E}P\bigl(\omega_0 \bigm| \omega_{-\infty}^{-1}\bigr)\;=\;1
\quad\forall\omega_{-\infty}^{-1}\in\Omega_{-\infty}^{-1}\, ,
\ee
such that:
\begin{itemize}
\item
the function $P(\omega_0 \,|\, \cdot\,)$ is \embf{continuous} for each 
$\omega_0\in E$, i.e., for all $\epsilon>0$ there exists $n\geq 0$ so that
\be{contlis}
\bigl|P\bigl(\omega_0 \bigm| \omega_{-\infty}^{-1}\bigr)
-P\bigl(\omega_0 \bigm| \sigma_{-\infty}^{-1}\bigr)\bigr| \;<\; \epsilon
\ee
for all $\omega_{-\infty}^{0},\sigma_{-\infty}^{0}\in\Omega_{-\infty}^{0}$ with
$\omega_{-n}^{0}=\sigma_{-n}^{0}$;
\item
the function $P(\omega_0 \,|\, \cdot\,)$ is \embf{strongly non-null} for each 
$\omega_0\in E$, i.e., $P(\omega_0 \,|\, \cdot\,)\geq c>0$.
\end{itemize}
\ed
[Property \reff{contlis} is indeed continuity with respect to the product discrete topology 
of $\Omega$.]

\bd{chain}
A probability measure $\mu$ on $(\Omega,\, \cF)$ with underlying process 
$(X_i)_{i\in\Z}$ on $(\Omega,\cF,\mu)$ 
is said to be a \embf{regular $g$-measure} if $\mu$ is shift-invariant and there exists 
a regular $g$-function $P$ such that $\mu$ is \embf{consistent} with $P$, namely,
\be{chainconsist}
\mu\bigl(X_{0}=\omega_{0} \bigm| X_{-\infty}^{-1}=\omega_{-\infty}^{-1}\bigr)
\;=\;P\bigl(\omega_{0} \bigm| \omega_{-\infty}^{-1}\bigr)
\ee
for all $\omega_{0}\in E$ and $\mu$-a.e.\ $\omega_{-\infty}^{-1}\in\Omega_{-\infty}^{-1}$. 
\ed
\br{extrem}
In the consistency definition (\ref{chainconsist}), $\mu$ needs only to be defined on
$(\Omega_{-\infty}^{0},\, \cF_{\leq 0})$. By shift-invariance, $\mu$ can be extended 
in a unique way to $(\Omega,\, \cF)$. Thus, without loss of generality, we identify 
$\mu$ on $(\Omega_{-\infty}^{0},\, \cF_{\leq 0})$ with its natural extension on 
$(\Omega,\, \cF)$.
\er

The previous definition is not very useful to prove the regularity of a measure.  For this, 
the following well known result is often useful.  Let us call a measure $\mu$ on 
$(\Omega,\, \cF)$ \embf{non-null} if it gives non-zero measure to every cylinder.  
We then have (see, for instance, \cite{kal90}):
\bt{th:2}
A probability measure $\mu$ on $(\Omega,\, \cF)$ with underlying process 
$(X_i)_{i\in\Z}$ on $(\Omega,\cF,\mu)$ 
is a regular $g$-measure iff it is non-null, shift-invariant and the sequences
\be{eq:r1}
\Bigl[\mu\bigl(X_{0}=\omega_{0} \bigm| X_{-n}^{-1}=\omega_{-n}^{-1}\bigr)\Bigr]_{n\ge 1} 
\quad\mbox{converge uniformly as } n\to \infty\;.
\ee
\et
\medskip

Let us now review the notions of specification and Gibbs measures.  It involves definitions 
and properties analogous to the preceding ones, but replacing one-side by two-side 
conditioning. (We specialize to the one-dimensional setting but similar notions and results 
are valid for any dimension.)
\bd{spe3}
A \emph{specification} $\gamma$ on $\left(\Omega,\mathcal{F}\right)$ is a family of 
probability kernels $\{\gamma_{\Lambda}\colon  \Lambda \subset\Z, |\Lambda|<\infty\}$, 
$\gamma_{\Lambda} \colon \mathcal{F} \times \Omega \rightarrow [0,1]$ such that for 
all finite $\Lambda\subset\Z$:
\begin{itemize}
\item
$\gamma_{\Lambda}(A \,|\, \cdot\,)\in \mathcal{F}_{\Lambda^{\cc}}$, 
for each $A \in \mathcal{F}$;
\item
$\gamma_{\Lambda}(B \,|\, \omega) = \one_{B}(\omega)$,
for each $B \in \mathcal{F}_{ \Lambda^{\cc}}$ and  $\omega \in \Omega$;
\item
$\gamma_{\Delta}\gamma_{\Lambda} \;=\; \gamma_{\Delta}$,
for each finite $\Delta \subset\Z\colon \Delta \supset \Lambda$, i.e.,
\be{gibbs5.1}
\iint h(\xi) \gamma_{\Lambda}(d \xi \,|\, \sigma) 
\gamma_{\Delta}(d \sigma \,|\, \omega) 
\;=\;\int h(\sigma)\gamma_{\Delta}(d \sigma \,|\, \omega)
\ee
for all measurable functions $h$ and configurations $\omega \in \Omega$.
\end{itemize}
\ed

In words, a specification is almost a \emph{regular system of  conditional probabilities} for 
finite regions conditioned on the external sigma algebras.  The only difference lies in the 
fact that conditional probabilities satisfy \reff{gibbs5.1} for \emph{almost all} $\omega$ with 
respect to some measure known beforehand.  In the definition of specification there is no 
previously known measure, thus condition  \reff{gibbs5.1} is asked \emph{for all} $\omega$
(see e.g.\ \cite{vEFS_JSP} for more details). We can now define regularity, as expected, 
as a property of all the kernels $\gamma_\Lambda$.  In our setting, however, we are 
interested only in non-null specifications and they, in fact, are entirely determined by their 
single-site part $\{\gamma_{\{i\}}\colon i\in\Z\}$ (see~\cite{fermai06} and references therein).  
Relevant definitions need only be done, thus, at the level of these single-site kernels.  That 
is what we do now, to emphasize the parallelism with the treatment of $g$-measures. For 
brevity we denote $\gamma_{\{i\}}(\omega_i \,|\, \cdot\,)$  the kernel $\gamma_{\{i\}}$ applied 
to the cylinder of base $\omega_i$.

\begin{definition}
\label{spedef}
A specification $\gamma$ on $(\Omega,\cF)$ is \embf{regular} if, for each $i\in\Z$,
\begin{itemize}
\item 
the function $\gamma_{\{i\}}(\omega_i \,|\, \cdot\,)$ is \embf{continuous} 
for each $\omega_i\in E$, i.e., for all $\epsilon>0$, there exists $n,m\geq 0$ so that
\be{contspe}
\bigl|\gamma_{\{i\}}\bigl(\omega_i \bigm| \omega_{\{i\}^{\cc}}\bigr)
-\gamma_{\{i\}}\bigl(\omega_i \bigm| \sigma_{\{i\}^{\cc}}\bigr)\bigr| \;<\;\epsilon
\ee
for all $\omega,\sigma\in\Omega$ with $\omega_{-n}^{m}=\sigma_{-n}^{m}$;
\item
the function $\gamma_{\{i\}}(\omega_i \,|\, \cdot\,)$ is \embf{strongly non-null} for 
each $\omega_i\in E$, i.e., $\gamma_{\{i\}}(\omega_i \,|\, \cdot\,)\geq c>0$.
\end{itemize}
\end{definition}

\begin{definition}
\label{d.con}
A probability measure $\mu$ on $(\Omega, \cF)$ with underlying field $(X_i)_{i\in\Z}$ 
on $(\Omega,\cF,\mu)$ is said to be a \embf{Gibbs measure} if $\mu$
is \embf{consistent} with some regular specification $\gamma$, namely,
\be{gibbsconsist}
\mu\bigl(X_i=\omega_i \bigm| X_{\{i\}^{\cc}}=\omega_{\{i\}^{\cc}}\bigr) \;=\;
\gamma_{\{i\}}\bigl(\omega_i \bigm| \omega_{\{i\}^{\cc}}\bigr)
\ee
for all $\omega_i\in E$ and $\mu$-a.e.\ $\omega_{\{i\}^{\cc}}\in\Omega_{\{i\}^{\cc}}$.
\end{definition}
The following is the analogues of Theorem \ref{th:2}.
\bt{th:3}
A shift-invariant probability measure $\mu$ on $(\Omega,\, \cF)$ with underlying field
$(X_i)_{i\in\Z}$ on $(\Omega,\cF,\mu)$ is a Gibbs measure iff it is non-null and the 
sequences
\be{eq:r3}
\Bigl[\mu\bigl(X_{0}=\omega_{0} \bigm| X_{-m}^{-1}
=\omega_{-m}^{-1}\,,\, X_1^m=\omega_1^m\bigr)\Bigr]_{n,m\ge 1} \quad
\mbox{converge uniformly as } n,m\to \infty\;.
\ee
\et
A measure violating \reff{eq:r3} is, thus, a \emph{non-Gibbsian measure} (see, e.g.\ \cite{fer06} 
and references therein).  More specifically, the absence of convergence of a sequence \reff{eq:r3} 
corresponds to an \emph{essential discontinuity} at $\omega\in\Omega$, that is, a discontinuity in 
the conditional expectations that can not be removed by a redefinition on a zero-measure set (a more detailed discussion of these issues can be found in \cite[Section 5.3]{fer06}). 

The link between Definition \ref{d.con} and the usual notion in classical statistical mechanics is 
provided by a theorem due to Kozlov~\cite{koz74} that states that a specification is Gibbsian if 
and only if it has the Boltzmann form 
\be{}
\gamma_\Lambda\bigl(\omega_i\bigm| \omega_{\{i\}^{\cc}}\bigr)
\;=\;\exp\Bigl\{-\sum_{A\ni i} \phi_A(\omega_A)\Bigr\}/\mbox{Norm.}\; ,
\ee
where the functions $\phi_A$ (\emph{interaction}) satisfy the summability condition
\be{eq:r-sum}
\sup_{i\in\Z}\sum_{A\ni i} \|\phi_A\|_\infty \;<\;\infty\;.
\ee
In the theory of dynamical systems, often Gibbsianness is associated with the SRB measures.  
These are measures $\mu$ for which there exists a function $\psi\colon \{-1,1\}^{\Z_+}\to\mathbb R$,
a constant $\Theta=\Theta(\psi)$, and some finite positive constants $\underline c,\overline c$ 
such that
\be{BowenGibbs}
\underline c\;\le\;
\frac {\mu(\omega_0^n)}{\exp\Bigl(\sum_{i=0}^{n}\psi(\tau^i \omega)-(n+1)\Theta\Bigr)}
\;\le\; \overline c\;
\ee
[$\tau^i$ is the $i$th iterate of the shift on $\Omega$]. In the non-null case, these SRB measures 
form a strict subset of the one-dimensional Gibbs measures.  This can be seen in two ways.  First, 
general Gibbs measures satisfy \reff{BowenGibbs} but with the constants $\underline c,\overline c$ 
substituted by $o(n)$-functions.  Second, the corresponding interactions for the non-null SRB 
measures must satisfy the condition
\be{eq:r-sum-srb}
\sup_{i\in\Z}\sum_{A\ni i} {\rm diam}(A)\,\|\phi_A\|_\infty \;<\;\infty
\ee
which is stronger than \reff{eq:r-sum}.


\section{Main result}
\label{S3}

For any 
$\underline{\omega}=\omega_{-\infty}^{-1}\in \Omega_{-\infty}^{-1}$
and $\overline{\omega}=\omega_{1}^{+\infty}\in \Omega_{1}^{+\infty}$, let 
\be{}
\ell(\underline{\omega}) \;=\; \min\{j\geq0\colon \omega_{-j-1}=1\}
\quad\text{and}\quad
m(\overline{\omega}) \;=\; \min\{j\geq0\colon \omega_{j+1}=1\}
\ee
denote the number of $0$'s before finding the first $1$ when looking backward in 
$\underline{\omega}$ and forward in $\overline{\omega}$, respectively 
[$\ell(\underline{0})=m(\overline{0})=\infty$].

Our $g$-function is defined by a converging sequence $\{p_{i}\}_{i\geq0}$ of numbers 
with values in $(0,1)$, satisfying
\be{eq:s1}
\inf_{i\geq0}p_{i}=\epsilon>0
\quad,\qquad p_{\infty}=\lim_{i\rightarrow+\infty}p_{i}\;.
\ee
The kernel $P$ is defined on $E\times\Omega_{-\infty}^{-1}$ by 
\be{eq:m1}
P(1\,|\, \underline{\omega}) \;=\; p_{\ell(\underline{\omega})}
\quad\forall \underline{\omega}\in\Omega_{-\infty}^{-1}.
\ee  
Note that the continuity of the kernel $P$ follows from the fact that
\be{}
\begin{aligned}
\sup_{{\underline{\omega},\underline{\sigma}\in\Omega_{-\infty}^{-1}}
\atop {\omega_{-k}^{-1}=\sigma_{-k}^{-1}}}
\Big|P\bigl(1 \bigm| \underline{\omega}\bigr)-P\bigl(1 \bigm| \underline{\sigma}\bigr)\Bigr|
&=\sup_{\omega_{-\infty}^{-k-1},\sigma_{-\infty}^{-k-1}\in\Omega_{-\infty}^{-k-1}}
\Bigl|P\bigl(1 \bigm| 0_{-k}^{-1}\omega_{-\infty}^{-k-1}\bigr)
-P\bigl(1 \bigm| 0_{-k}^{-1}\sigma_{-\infty}^{-k-1}\bigr)\Bigr|\\
&=\sup_{l,m \geq k}|p_{l}-p_{m}|
\end{aligned}
\ee
which goes to $0$ since $\{p_{k}\}_{k\geq0}$ converges.

This $g$-function is, therefore, regular and, furthermore, in~\cite{gal11} it is proven
that there exists a unique stationary chain $\mu$ compatible with $P$ which is
the renewal chain with infinitely many $1$'s separated by intervals of $0$'s 
of random length and having exponential tail distribution.   For all practical purposes, 
this chain is as regular as it can be.  Nevertheless, it is not necessary Gibbsian.

\bt{}\label{th:1}
There exist choices of the sequence $\{p_{i}\}_{i\geq0}$ satisfying {\rm\reff{eq:s1}} for which 
\be{eq:rr09}
\mu\bigl(X_0=0 \bigm|\cdot\,\bigr) \mbox{ is essentially discontinuous at } 0_{-\infty}^{+\infty}\;,
\ee
where $\mu$ is the (unique and non-null) $g$-measure compatible with the kernel defined by 
{\rm\reff{eq:m1}}. 
\et

By Kozlov's theorem \cite{koz74} this means that the resulting $g$-measure is non-Gibbsian in 
the statistical mechanical sense, and hence neither in SRB sense.

\bpr
The proof consists in the observation that if $\omega$ has $\ell(\underline{\omega})=i$ and $m(\overline{\omega})=j$, then, as we will see in \reff{ratioeq} below, 
$\mu(X_{0}=0 \,|\, X_{-m}^{-1}=\omega_{-m}^{-1}\,,\, X_{1}^{n}=\omega_{1}^{n})$
is determined by the ratio
\be{}
\prod_{k=0}^{i-1}\frac{1-p_{k}}{1-p_{k+j}}.
\ee
Thus, the discontinuity at $0_{-\infty}^{+\infty}$ is equivalent to the existence of 
a sequence of $p_k$ for which this ratio oscillates with $i$ and $j$.  The most 
economical way of achieving this is to define
\be{pkdef}
p_{k} \;=\; 1-(1-p_{\infty})\xi^{v_{k}}\; ,
\ee 
so that 
\be{eq:rr10}
\prod_{k=0}^{i-1}\frac{1-p_{k}}{1-p_{k+j}} \;=\; \xi^{\sum_{k=0}^{i-1} (v_{k}-v_{k+j})}\;.
\ee
The discontinuity is obtained by choosing a sequence  $v_{k}$ converging to $0$ when $k\to\infty$, but such that 
$\sum_{k=0}^{i}v_{k}$ oscillates.
\medskip

To formalize this idea, let us first provide an explicit expression for the conditional probabilities.
By construction, the measure $\mu$ compatible with $P$ has the property that, if
$\ell(\underline{\omega})=i\leq m<\infty$ and $m(\overline{\omega})=j\leq n<\infty$ 
\be{}
\begin{aligned}
&\mu\bigl(X_{-m}^{n}=\omega_{-m}^{n}\bigr)\\
&\quad=\; \mu\bigl(X_{-m}^{-i-1}=\omega_{-m}^{-i-1}\bigr)
\mu\bigl(X_{-i}^{j+1}=\omega_{-i}^{j+1} \bigm| X_{-i-1}=1\bigr)
\mu\bigl(X_{j+2}^{n}=\omega_{j+2}^{n} \bigm| X_{j+1}=1\bigr)\; 
\end{aligned}
\ee
which becomes 
\begin{equation}
\mu\bigl(X_{-m}^{-i-1}=\omega_{-m}^{-i-1}\bigr)
\biggl(\prod_{k=0}^{i+j}(1-p_{k})p_{i+j+1}\biggr)
\mu\bigl(X_{j+2}^{n}=\omega_{j+2}^{n} \bigm| X_{j+1}=1\bigr)
\end{equation}
when $\omega_{0}=0$, and
\begin{equation}
\mu\bigl(X_{-m}^{-i-1}=\omega_{-m}^{-i-1}\bigr)
\biggl(\prod_{k=0}^{i-1}(1-p_{k})p_{i}\biggr)
\biggl(\prod_{k=0}^{j-1}(1-p_{k})p_{j}\biggr)
\mu\bigl(X_{j+2}^{n}=\omega_{j+2}^{n} \bigm| X_{j+1}=1\bigr)
\end{equation}
when $\omega_{0}=1$.
Thus, the finite-volume 2-sided conditional probability equals
\be{ratioeq}
\begin{aligned}
&\mu\bigl(X_{0}=0 \bigm| X_{-m}^{-1}=\omega_{-m}^{-1},\, X_{1}^{n}=\omega_{1}^{n}\bigr)\\
&\quad=\; \frac{\mu\bigl(X_{-m}^{n}=\omega_{-m}^{n}\bigr)}
{\mu\bigl(X_{-m}^{n}=\omega_{-m}^{n}\bigr)
+\mu\bigl(X_{-m}^{n}=\omega_{-m}^{-1}1_{0}\omega_{1}^{n}\bigr)}\\
&\quad=\; \frac{\prod_{k=0}^{i+j}(1-p_{k})p_{i+j+1}}
{\prod_{k=0}^{i-1}(1-p_{k})p_{i}
\prod_{k=0}^{j-1}(1-p_{k})p_{j}+\prod_{k=0}^{i+j}(1-p_{k})p_{i+j+1}}\\
&\quad=\; \Biggl(1+\frac{p_{i}p_{j}}{(1-p_{i+j})p_{i+j+1}} 
\prod_{k=0}^{i-1}\frac{1-p_{k}}{1-p_{k+j}}\Biggr)^{-1}
\end{aligned}
\ee
for all $n\ge i$ and $m\ge j$.  Notice that the intermediate equality is symmetrtic under the interchange $i\leftrightarrow j$; this proves the left-right symmetry of the conditional expectations.
Writing the probabilities $p_k$ in the form \reff{pkdef} this becomes
\be{eq:mu1}
\begin{aligned}
&\mu\bigl(X_{0}=0 \bigm| X_{-i-1}^{-1}=1_{-i-1}0_{-i}^{-1},\, 
X_{1}^{j+1}=0_{1}^{j}1_{j+1}\bigr)\\
&\qquad=\; \Biggl(1+\frac{[1-(1-p_{\infty})\xi^{v_{i}})][1-(1-p_{\infty})\xi^{v_{j}})]}
{(1-p_{\infty})\xi^{v_{i+j}}[1-(1-p_{\infty})\xi^{v_{i+j+1}})]} 
\xi^{\sum_{k=0}^{i-1} (v_{k}-v_{k+j})}\Biggr)^{-1}\; .
\end{aligned}
\ee
\medskip

To conclude we must choose a sequence $v_k$ with an oscillating sum, as proposed after 
\reff{eq:rr10}.  We choose $\xi\in(1,(1-p_{\infty})^{-2})$ and consider the sequence $\{v_{k}\}_{k\ge0}$ with
\be{vk:def}
v_{k} \;=\; \frac{(-1)^{r_{k}}}{r_{k}}
\quad\text{with}\quad 
r_{k} \;=\; \inf\biggl\{i\geq1\colon\sum_{j=1}^{i}j\geq k+1\biggr\}\; .
\ee 
The first terms of this sequence are as follows:
\be{}
-1\,,\,\,\frac{1}{2}\,,\,\,\frac{1}{2}\,,\,\,-\frac{1}{3}\,,\,\,-\frac{1}{3}\,,\,\,-\frac{1}{3}\,,\,\,
\frac{1}{4}\,,\,\,\frac{1}{4}\,,\,\,\frac{1}{4}\,,\,\,\frac{1}{4}\,,\,\,-\frac{1}{5}\,,\,\,
-\frac{1}{5}\,,\,\,-\frac{1}{5}\,,\,\,-\frac{1}{5}\,,\,\,-\frac{1}{5}\,,\,\,\ldots
\ee
Clearly, $v_{k}$ converges to $0$ when $k\to\infty$, while $\sum_{k=0}^{i}v_{k}$ 
does not converge as $i$ diverges because it oscillates inside $[-1,0]$.  
Due to \reff{vk:def}, there are strictly increasing subsequences
\begin{itemize}
\item $\{i_{n}^{(1)}\}_{n\geq 1}$ such that $\sum_{k=0}^{i_{n}^{(1)}-1}v_{k}=-1$ 
for any $n\geq 1$;
\item $\{j_{n}^{(1)}\}_{n\geq 1}$, where $j_{n}^{(1)}$ depends on $i_{n}^{(1)}$, such that
$\sum_{k=0}^{i_{n}^{(1)}-1}v_{j_{n}^{(1)}+k}=0$ for any $n\geq 1$;
\item $\{i_{n}^{(2)}\}_{n\geq 1}$ such that $\sum_{k=0}^{i_{n}^{(2)}-1}v_{k}=0$ 
for any $n\geq 1$;
\item $\{j_{n}^{(2)}\}_{n\geq 1}$, where $j_{n}^{(2)}$ depends on $i_{n}^{(2)}$, such that
$\sum_{k=0}^{i_{n}^{(2)}-1}v_{j_{n}^{(2)}+k}=0$ for any $n\geq 1$.
\end{itemize}
Therefore, we have that
\be{}
\sum_{k=0}^{i_{n}^{(1)}-1}\Bigl(v_{k}-v_{j_{n}^{(1)}+k}\Bigr) \;=\; -1
\ee
and
\be{}
\sum_{k=0}^{i_{n}^{(2)}-1}\Bigl(v_{k}-v_{j_{n}^{(2)}+k}\Bigr) \;=\; 0\; .
\ee
Since
\be{eq:mu5}
\lim_{i,j\to\infty} \frac{[1-(1-p_{\infty})\xi^{v_{i}})][1-(1-p_{\infty})\xi^{v_{j}})]}
{(1-p_{\infty})\xi^{v_{i+j}}[1-(1-p_{\infty})\xi^{v_{i+j+1}})]} 
\;=\; \frac{p_{\infty}}{(1-p_{\infty})}\; ,
\ee
it follows from \reff{eq:mu1}--\reff{eq:mu5} that
\begin{eqnarray}
\lefteqn{\hspace{-3cm}\lim_{n\to\infty}\mu\Bigl(X_{0}=0 \Bigm | X_{-i^{(\ell)}_n-1}^{-1}=1_{-i^{(\ell)}_n-1}
0_{-i^{(\ell)}_n}^{-1},\, 
X_{1}^{j^{(\ell)}_n+1}=0_{1}^{j^{(\ell)}_n}1_{j^{(\ell)}_n+1}\Bigr)}\nonumber\\[10pt]
&&=\;\left\{\begin{array}{ll}
\displaystyle\frac{(1-p_\infty)\xi}{(1-p_\infty)\xi+p_\infty} & \mbox{if } \ell=1\\[10pt]
(1-p_\infty) & \mbox{if } \ell=2\;.
\end{array}\right.
\end{eqnarray}
Hence $\lim_{i,j\to\infty}\mu\bigl(X_{0}=0 \bigm| X_{-i-1}^{-1}=1_{-i-1}0_{-i}^{-1},\, 
X_{1}^{j+1}=0_{1}^{j}1_{j+1}\bigr)$
does not exist.  This contradicts \reff{eq:r3} and, more specifically, proves the essential discontinuity \reff{eq:rr09}.
\epr




\begin{thebibliography}{99}



\bibitem{bow75}
R.~Bowen.
{\it Equilibrium states and the ergodic theory of {A}nosov diffeomorphisms},
Lecture Notes in Mathematics, Vol.\ 470, Springer-Verlag, Berlin, 1975.

\bibitem{dob68}
R.~L.\ Dobrushin,
The description of a random field by means of conditional probabilities and 
conditions of its regularity.
{\it Theory Probab.\ Appl.}, {\bf 13} (1968) 197--224.

\bibitem{vEFS_JSP}
A.C.D.\ van Enter, R.\ Fern{\'a}ndez, and A.D.\ Sokal,
Regularity properties and pathologies of position-space renormalization-group 
transformations: scope and limitations of {G}ibbsian theory,
{\it J.\ Stat.\ Phys.} {\bf 72} (1993) 879--1167.

\bibitem{entver10}
A.C.D.\ van Enter and E.\ Verbitskiy, 
Erasure entropies and Gibbs measures,
{\it Markov Process.\ Related Fields} {\bf 16} (2010) 3--14.

\bibitem{fer06}
R.\ Fern\'andez, 
Gibbsianness and non-Gibbsianness in Lattice random fields, 
In: {\it Mathematical statistical physics}, (A.\ Bovier, A.C.D.\ van Enter, F.\ den Hollander
and F.\ Dunlop eds.), 731--799, Elsevier B.\ V., Amsterdam, 2006.

\bibitem{fermai04}
R.\ Fern\'andez and G.\ Maillard,
Chains with complete connections and one-dimensional Gibbs measures, 
{\it Electron.\ J.\ Probab.} {\bf 9} (2004) 145--176.



\bibitem{fermai06}
R.\ Fern\'andez and G.\ Maillard,
Construction of a specification from its singleton part,
{\it ALEA} {\bf 2} (2006) 297--315.

\bibitem{gal11}
S.\ Gallo,
Chains with unbounded variable length memory:
perfect simulation and visible regeneration scheme,
{\it To appear in Adv.\ in Appl.\  Probab.}

\bibitem{geo88}
H.-O. Georgii,
{\it Gibbs Measures and Phase Transitions},
de Gruyter Studies in Mathematics, Vol.\ 9, 
Walter de Gruyter \& Co., Berlin, 1988.

\bibitem{har55}
T.~E.\ Harris,
On chains of infinite order,
{\it Pacific J.\ Math.}, {\bf 5} (1955) 707--24.

\bibitem{kal90}
S.\ Kalikow,
Random {M}arkov processes and uniform martingales,
{\it Isr.\ J.\ Math.}, {\bf 71} (1990) 33--54.

\bibitem{kea72}
M.\ Keane,
Strongly mixing $g$-measures,
{\it Inventiones Math.} {\bf 16} (1972) 309--324.

\bibitem{kel98}
G.\ Keller,
{\it Equilibrium states in ergodic theory}, 
London  Mathematical Society Student Texts, Vol.\ 42,
Cambridge University Press, Cambridge, 1998.

\bibitem{koz74}
O.~K.\ Kozlov,
Gibbs description of a system of random variables,
{\it Probl.\ Inform.\ Transmission}, {\bf 10} (1974) 258--265.

\bibitem{lanrue69}
O.~E.\ {Lanford III} and D.\ Ruelle,
Observables at infinity and states with short range correlations in statistical mechanics,
{\it Comm.\ Math.\ Phys.}, {\bf 13} (1969) 194--215.


\bibitem{onimih35}
O.\ Onicescu and G.\ Mihoc,
Sur les cha\^{\i}nes statistiques,
{\it C.\ R.\ Acad. Sci.\ Paris}, {\bf 200} (1935) 511--512.

\bibitem{rue78}
D.\ Ruelle,
{\it Thermodynamic formalism}.
Encyclopedia of {M}athematics {\bf 5}, Addison Wesley, New-York, 1978.


\end{thebibliography}
\end{document}